\documentclass{amsart}
\usepackage{amsmath}
\usepackage{amstext}
\usepackage{amssymb}
\usepackage{amsthm}
\usepackage{amscd}

\usepackage{upref}

\swapnumbers
\theoremstyle{plain}
\newtheorem{thm}{Theorem}[section]
\newtheorem{lem}[thm]{Lemma}
\newtheorem{prop}[thm]{Proposition}
\newtheorem{cor}[thm]{Corollary}

\theoremstyle{definition}
\newtheorem{defi}[thm]{Definition}

\theoremstyle{remark}

\newtheorem{rmk}[thm]{Remark}

\DeclareMathOperator{\Coker}{Coker} \DeclareMathOperator{\Ima}{Im}

\DeclareMathOperator{\rank}{rank}

\DeclareMathOperator{\length}{length}
\DeclareMathOperator{\content}{c}

\def\Z{\mathbb Z}
\def\N{\mathbb N}
\def\Q{\mathbb Q}

\def\nn{\relax\ifmmode{\mathbb N_{0}}\else$\mathbb N_{0}$\fi}
\def\lra{\longrightarrow}

\begin{document}

\title[Some properties of top graded local cohomology modules]
{Some properties of top graded local cohomology modules}
\author{MORDECHAI KATZMAN}
\address[Katzman]{Department of Pure Mathematics,
University of Sheffield, Hicks Building, Sheffield S3 7RH, United Kingdom\\
{\it Fax number}: 0044-114-222-3769}
\email{M.Katzman@sheffield.ac.uk}
\author{RODNEY Y. SHARP}
\address[Sharp]{Department of Pure Mathematics,
University of Sheffield, Hicks Building, Sheffield S3 7RH, United Kingdom\\
{\it Fax number}: 0044-114-222-3769}
\email{R.Y.Sharp@sheffield.ac.uk}

\subjclass{Primary 13D45, 13E05, 13A02, 13P10; Secondary 13D40}

\date{\today}

\keywords{Graded commutative Noetherian ring, graded local
cohomology module, module of inverse polynomials, associated prime
ideal, cohomological Hilbert function}


\maketitle

\setcounter{section}{-1}
\section{\bf Introduction}
\label{in}

Let $R = \bigoplus_{d \in \nn} R_d$ be a positively graded
commutative Noetherian ring which is standard in the sense that $R
= R_0[R_1]$, and set $R_+ := \bigoplus_{d \in \N} R_d$, the
irrelevant ideal of $R$. (Here, $\nn$ and $\N$ denote the set of
non-negative and positive integers respectively; $\Z$ will denote
the set of all integers.) Let $M = \bigoplus_{d \in \Z} M_d$ be a
non-zero finitely generated graded $R$-module. This paper is
concerned with the behaviour of the graded components of the
graded local cohomology modules $H^i_{R_+}(M)~(i \in \nn)$ of $M$
with respect to $R_+$.

It is known (see \cite[15.1.5]{BS}) that there exists $r \in \Z$
such that $H^i_{R_+}(M)_{d} = 0$ for all $i \in \nn$ and all $d
\geq r$, and that $H^i_{R_+}(M)_{d}$ is a finitely generated
$R_0$-module for all $i \in \nn$ and all $d \in \Z$.

The first part (\S\ref{section1}) of this paper deals with the
case in which $R=R_0[U_1, \dots, U_s]/I$, where $U_1, \dots, U_s$
are indeterminates of degree one, and $I\subset R_0[U_1, \dots,
U_s]$ is a homogeneous ideal. The main theorem of that section is
that for $d\geq s$, all the associated primes of
$H^s_{R_+}(R)_{-d}$ contain a certain ideal of $R_0$ called the
``content'' of $I$ (see Definition \ref{definition1}.) This result
provides an affirmative answer, in a special case, to a question
discussed by M. Brodmann and M. Hellus in \cite{BH}, namely,
whether it is true that there exists $d_0 \in \Z$ such that either
$H^i_{R_+}(M)_{-d}=0$ for all $d>d_0$ or $H^i_{R_+}(M)_{-d}\neq 0$
for all $d>d_0$. In \cite[Lemma (4.2)]{BH}, Brodmann and Hellus
showed that the answer is affirmative for all values of $i$ when
$R_0$ is semi-local and of dimension not exceeding $1$; in \S1
below, we show that the answer is affirmative when $R_0$ is a
domain, $M = R$ and $R_+$ can be generated by $i$ homogeneous
elements of degree $1$.

There are instances when all the homogeneous components
$H^i_{R_+}(M)_{r}~(r \in \Z)$ have finite length (for example,
when $R_0$ is Artinian) and we may then define the \emph{$i$-th
cohomological Hilbert function of M,} denoted by $h_M^i : \Z \lra
\nn$, by $h_M^i(r)= \length _{R_0}H^i_{R_+}(M)_{r}$ for all $r \in
\Z$. When $R_0$ is Artinian this function agrees with a polynomial
for all $r < < 0$ (see \cite[Theorem 17.1.9]{BS}). In \S\S2,3 of
this paper we construct examples which show that this result need
not be true when $R_0$ is not Artinian.

\section{\bf The vanishing of top local cohomology modules}\label{section1}

Throughout this section $R_0$ will denote an arbitrary commutative
Noetherian domain. We set $S=R_0[U_1, \dots, U_s]$ where $U_1,
\dots, U_s$ are indeterminates of degree one, and $R=S/I$ where
$I\subset R_0[U_1, \dots, U_s]$ is a graded ideal. For $t \in
\Z$, we shall denote by $(\: {\scriptscriptstyle \bullet} \:)(t)$
the $t$-th shift functor (on the category of graded $R$-modules
and homogeneous homomorphisms).

For any multi-index $\lambda=(\lambda^{(1)}, \dots, \lambda^{(s)})
\in \Z^s$ we shall write $U^\lambda$ for $U_1^{\lambda^{(1)}}
\dots U_s^{\lambda^{(s)}}$ and we shall set $| \lambda | =
\lambda^{(1)}+ \dots+ {\lambda^{(s)}}$.

\begin{lem}\label{lemma1}
Let $I$ be generated by homogeneous elements $f_1, \dots, f_r\in
S$. Then there is an exact sequence of graded $S$-modules and
homogeneous homomorphisms
$$ \bigoplus_{i=1}^r H^s_{S_+}(S)(-\deg f_i)\xrightarrow[]{(f_1, \dots, f_r)}
H^s_{S_+}(S) \longrightarrow H^s_{R_+}(R) \longrightarrow 0.$$
\end{lem}

\begin{proof}
The functor $H^s_{S_+}$ is right exact and the natural equivalence
between $H^s_{S_+}$ and $(\: {\scriptscriptstyle \bullet}
\:)\otimes_S H^s_{S_+}(S)$ (see \cite[6.1.8 \& 6.1.9]{BS})
actually yields a homogeneous $S$-isomorphism
$$H^s_{S_+}(S)/(f_1,\ldots,f_r)H^s_{S_+}(S) \cong H^s_{S_+}(R).$$ To complete the
proof, just note that there is an isomorphism of graded
$S$-modules $H^s_{S_+}(R) \cong H^s_{R_+}(R)$, by the Graded
Independence Theorem \cite[13.1.6]{BS}.
\end{proof}

We can realize $H^s_{S_+}(S)$ as the module $R_0[U_1^-,\dots,
U_s^-]$ of inverse polynomials described in \cite[12.4.1]{BS}:
this graded $R$-module has end $-s$, and, for each $d \geq s$, its
$(-d)$-th component is a free $R_0$-module with base $
{\mathcal{B}(d)}:=\left(U^\lambda\right)_{-\lambda\in \N^s,
|\lambda| = -d}. $ We combine this realisation with the previous
lemma to find a presentation of each homogeneous component of
$H^s_{R_+}(R)$ as a cokernel of a matrix with entries in $R_0$.

Assume first that $I$ is generated by one homogeneous element $f$
of degree $\delta$. For any $d$ we have, in view of Lemma
\ref{lemma1}, a graded exact sequence
$$R_0[U_1^-,\dots, U_s^-]_{-d-\delta} \xrightarrow[]{\phi_d}
R_0[U_1^-,\dots, U_s^-]_{-d} \longrightarrow H^s_{R_+}(R)_{-d}
\longrightarrow 0 .$$ The map of free $R_0$-modules $\phi_d$ is given
by multiplication on the left by a ${{d-1} \choose s-1} \times
{{d+\delta-1} \choose {s-1}}$ matrix which we shall denote later
by $M(f;d)$.

In the general case, where $I$ is generated by homogeneous
elements $f_1, \dots, f_r\in S$, it follows from Lemma
\ref{lemma1} that the $R_0$-module $H^s_{R_+}(R)_{-d}$ is the
cokernel of a matrix $M(f_1, \dots, f_r; d)$ whose columns consist
of all the columns of $M(f_1,d), \dots, M(f_r,d)$.

Consider a homogeneous $f\in S$ of degree $\delta$. We shall now
describe the matrix $M(f;d)$ in more detail and to do so we start
by ordering the bases of the source and target of $\phi_d$ as
follows. For any $\lambda, \mu\in \Z^s$ with negative entries we
declare that $U^\lambda < U^\mu$ if and only if $U^{-\lambda}
<_{\mathrm{Lex}} U^{-\mu}$ where ``$<_{\mathrm{Lex}}$'' is the
lexicographical term ordering in $S$ with $U_1 > \dots > U_s$. We
order these bases, and by doing so also the columns and rows of
$M(f;d)$, in ascending order.

\begin{lem}\label{lemma2}
Let $f\neq 0$ be a homogeneous element in $S$. Then, for all
$d\geq s$, the matrix $M(f;d)$ has maximal rank.
\end{lem}
\begin{proof}
We prove the lemma by producing a non-zero maximal minor of
$M(f;d)$. Write $f=\sum_{\lambda\in \Lambda} a_\lambda U^\lambda$
where $a_\lambda\in R_0\setminus\{ 0 \}$ for all $\lambda\in
\Lambda$ and let $\lambda_0$ be such that $U^{\lambda_0}$ is the
minimal member of $\left\{U^{\lambda} : \lambda \in
\Lambda\right\}$ with respect to the lexicographical term order in $S$.

Let $\delta$ be the degree of $f$. Each column of $M(f;d)$
corresponds to a monomial $U^\lambda\in {\mathcal{B}(d+\delta)}$;
its $\rho$-th entry is the coefficient of  $U^\rho$ in $f
U^\lambda \in R_0[U_1^-,\dots, U_s^-]_{-d}$.

Fix any $U^{\nu}\in {\mathcal{B}(d)}$ and consider the column
$c_\nu$ corresponding to $U^{\nu-\lambda_0 }\in
{\mathcal{B}(d+\delta)}$. The $\nu$-th entry of $c_\nu$ is
obviously $a_{\lambda_0}$. Also, for any other $\lambda_1 \in
\Lambda$ with $U^{\lambda_1} >_{\mathrm{Lex}} U^{\lambda_0}$,
either $\nu-\lambda_0+\lambda_1$ has an entry in $\nn$, in which
case the corresponding term of $f U^{\nu-\lambda_0}\in
R_0[U_1^-,\dots, U_s^-]_{-d}$ is zero, or
$U^{\nu-\lambda_0+\lambda_1} < U^\nu$. This last statement follows
from the fact that if $j$ is the first coordinate where
$\lambda_0$ and $\lambda_1$ differ then
$\lambda_0^{(j)}<\lambda_1^{(j)}$ and so also $-\nu^{(j)}>
-\nu^{(j)}+\lambda_0^{(j)}-\lambda_1^{(j)}$; this implies that
$U^{-\nu+\lambda_0-\lambda_1} <_{\mathrm{Lex}} U^{-\nu}$ and
$U^{\nu-\lambda_0+\lambda_1} < U^\nu$. We have shown that the last
non-zero entry in $c_\nu$ occurs at the $\nu$-th row and is equal
to $a_{\lambda_0}$. Consider the square submatrix of $M(f;d)$
whose columns are the $c_\nu~(\nu\in {\mathcal{B}(d)})$; its
determinant is clearly a power of $a_{\lambda_0}$ and hence is
non-zero.
\end{proof}

\begin{defi}\label{definition1}
For any $f \in R_0[U_1, \dots, U_s]$ write $f=\sum_{\lambda\in
\Lambda} a_\lambda U^\lambda$ where $a_\lambda\in R_0$ for all
$\lambda\in \Lambda$. For such an $f \in R_0[U_1, \dots, U_s]$ we
define the {\em content\/} $\content(f)$ of $f$ to be the ideal
$\langle a_\lambda : \lambda\in \Lambda \rangle$ of $R_0$
generated by all the coefficients of $f$. If $J\subset R_0[U_1,
\dots, U_s]$ is an ideal, we define its {\em content\/}
$\content(J)$ to be the ideal of $R_0$ generated by the contents
of all the elements of $J$. It is easy to see that if $J$ is
generated by $f_1, \dots, f_r$, then $\content(J) =
\content(f_1)+\dots +\content(f_r)$.
\end{defi}

\begin{lem}\label{lemma3}
Suppose that $I$ is generated by homogeneous elements $f_1, \dots,
f_r\in S$. Fix any $d\geq s$. Let $t:=\rank M(f_1, \dots, f_r; d)
={{d-1} \choose s-1}$ and let $I_d$ be the ideal generated by all
$t\times t$ minors of $M(f_1, \dots, f_r; d)$. Then $\content(I)
\subseteq \sqrt{I_d}$.
\end{lem}
\begin{proof}
It is enough to prove the lemma when $r=1$; let $f=f_1$. Write
$f=\sum_{\lambda\in \Lambda} a_\lambda U^\lambda$ where
$a_\lambda\in R_0\setminus\{ 0 \}$ for all $\lambda\in \Lambda$.
Assume that $\content(I) \not\subseteq \sqrt{I_d}$ and pick
$\lambda_0$ so that $U^{\lambda_0}$ is the minimal element in
$\left\{U^{\lambda} : \lambda \in \Lambda\right\}$ (with respect
to the lexicographical term order in $S$) for which $a_\lambda \notin
\sqrt{I_d}$. Notice that the proof of Lemma \ref{lemma2} shows
that $U^{\lambda_0}$ cannot be the minimal element of
$\left\{U^{\lambda} : \lambda \in \Lambda\right\}$.

Fix any $U^{\nu}\in {\mathcal{B}(d)}$ and consider the column
$c_\nu$ corresponding to $U^{\nu-\lambda_0} \in
{\mathcal{B}(d+\delta)}$. The $\nu$-th entry of $c_\nu$ is
obviously $a_{\lambda_0}$. An argument similar to the one in the
proof of Lemma \ref{lemma2} shows that, for any other $\lambda_1
\in \Lambda$ with $U^{\lambda_1} >_{\mathrm{Lex}} U^{\lambda_0}$,
either $\nu-\lambda_0+\lambda_1$ has an entry in $\nn$, in which
case the corresponding term of $f U^{\nu-\lambda_0}\in
R_0[U_1^-,\dots, U_s^-]_{-d}$ is zero, or $U^\nu >
U^{\nu-\lambda_0+\lambda_1}$.

Similarly, for any other $\lambda_1 \in \Lambda$ with
$U^{\lambda_1} <_{\mathrm{Lex}} U^{\lambda_0}$, either
$\nu-\lambda_0+\lambda_1$ has an entry in $\nn$, in which case the
corresponding term of $f U^{\nu-\lambda_0}\in R_0[U_1^-,\dots,
U_s^-]_{-d}$ is zero, or $U^\nu < U^{\nu-\lambda_0+\lambda_1}$.

We have shown that all the entries above the $\nu$-th row of
$c_\nu$ are in $\sqrt{I_d}$. Consider the matrix $M$ whose columns
are $c_\nu~(\nu \in {\mathcal{B}(d)})$ and let
$\overline{\phantom{X} } : R_0 \rightarrow R_0 / \sqrt{I_d}$
denote the quotient map. We have
$$0=\overline{\det(M)}=\det(\overline{M})=\overline{a_{\lambda_0}}^{{d-1} \choose s-1}$$
and, therefore, $a_{\lambda_0} \in \sqrt{I_d}$, a contradiction.
\end{proof}

\begin{thm}\label{theorem1}
Suppose that $I$ is generated by homogeneous elements $f_1, \dots,
f_r\in S$. Fix any $d\geq s$. Then each associated prime of
$H^s_{R_+}(R)_{-d}$
contains $\content(I)$. In particular
$H^s_{R_+}(R)_{-d}=0$ if and only if $\content(I)=R_0$.
\end{thm}

\begin{proof}
Recall that for any $p,q \in \N$ with $p\leq q$ and any $p\times
q$ matrix $M$ of maximal rank with entries in any domain, $\Coker
M = 0$ if and only if the ideal generated by the maximal minors of
$M$ is the unit ideal. Let $M=M(f_1, \dots, f_r; d)$, so that
$H^s_{R_+}(R)_{-d}\cong\Coker M$.

In view of Lemmas \ref{lemma2} and \ref{lemma3}, the ideal
$\content(I)$ is contained in the radical of the ideal generated
by the maximal minors of $M$; therefore, for each $x\in
\content(I)$, the localization of $\Coker M$ at $x$ is zero;
we deduce that
$\content(I)$ is contained in all associated primes of $\Coker M$.

To prove the second statement, assume first that $\content(I)$ is
not the unit ideal. Since all minors of $M$ are contained in
$\content(I)$, these cannot generate the unit ideal and $\Coker M
\neq 0$. If, on the other hand,  $\content(I)=R_0$ then $\Coker M$
has no associated prime and $\Coker M = 0$.
\end{proof}

\begin{cor}\label{cor0} Let the situation be as in\/ {\rm \ref{theorem1}}. The
following statements are equivalent:
\begin{enumerate}
\item $\content(I)=R_0$;
\item $H^s_{R_+}(R)_{-d}=0$ for some $d \geq s$;
\item $H^s_{R_+}(R)_{-d}=0$ for all $d \geq s$.
\end{enumerate}
Consequently, $H^s_{R_+}(R)$ is asymptotically gap-free in the sense of\/
{\rm \cite[(4.1)]{BH}}.

It follows that, if $T = \bigoplus_{n \in \nn}T_n$ is any standard
graded finitely generated $R_0$-algebra with $T_0 = R_0$ and if
the $R_0$-module $T_1$ can be generated by $s$ elements, then
$H^s_{T_+}(T)$ is asymptotically gap-free, because there is a
homogeneous surjective ring homomorphism $S \rightarrow T$.
\end{cor}

\begin{rmk}
Theorem \ref{theorem1} cannot be extended to $H^i_{R_+}(R)_{-d}=0$
for all $i<s$. For example, take $R_0=\mathbb{C}[X_1, X_2]$,
$S=R_0[U_1, U_2]$ (where $X_1,X_2,U_1,U_2$ are indeterminates) and
$$I=\left ( U_1(X_1 U_1 + X_2 U_2), U_2(X_1 U_1 + X_2 U_2) \right)
\subseteq S.$$ Then $\content(I)=( X_1, X_2)$ but
$H^0_{( U_1, U_2 )} (S/I) \cong (X_1 U_1 + X_2 U_2 )/ I$ and it
does not vanish when localized at $X_1$.
\end{rmk}

We consider the following consequence of
our work to be interesting because of
its relevance to associated primes of local cohomology modules.

\begin{cor}\label{cor2}
The $R$-module $H^s_{R_+}(R)$ has finitely many minimal associated
primes, and these are just the minimal primes of the ideal
$\content(I)R + R_+$.
\end{cor}

\begin{proof}
Let $r \in \content(I)$. By Theorem \ref{theorem1}, the
localization of $H^s_{R_+}(R)$ at $r$ is zero. Hence each
associated prime of $H^s_{R_+}(R)$ contains $\content(I)R$.  Such
an associated prime must contain $R_+$, since $H^s_{R_+}(R)$ is
$R_+$-torsion.

On the other hand, $H^s_{R_+}(R)_{-s} \cong R_0/\content(I)$ and
$H^s_{R_+}(R)_i = 0$ for all $i > -s$; therefore there is an
element of the $(-s)$-th component of $H^s_{R_+}(R)$ that has
annihilator (over $R$) equal to $\content(I)R + R_+$. All the
claims now follow from these observations.
\end{proof}

\begin{rmk}\label{rmk1}
In \cite[Conjecture 5.1]{Hunek92}, C. Huneke put forward the
conjecture that every local cohomology module (with respect to any
ideal) of a finitely generated module over a local Noetherian ring
has only finitely many associated primes. Recently, this
conjecture was settled by the first author, who gave a
counterexample in \cite[Corollary 1.3]{Katzm01}.  Corollary
\ref{cor2} provides a little evidence in support of the weaker
conjecture that every local cohomology module (with respect to any
ideal) of a finitely generated module over a local Noetherian ring
has only finitely many {\em minimal\/} associated primes.
\end{rmk}

\section{\bf Cohomological Hilbert functions which are not of
reverse polynomial type}\label{sec2}

In this section, we shall use the terminology used in the final
paragraph of the Introduction; also, we follow the terminology of
\cite[17.1.1(vi)]{BS}, and say that a function $f : \Z \rightarrow
\Z$ is of {\it reverse polynomial type\/} if and only if there
exists a polynomial $P \in \Q[T]$ such that $f(r) = P(r)$ for all
$r < < 0$. In this and the next section, we are interested in
situations where $H^i_{R_+}(M)_{-d}$ has finite length as an
$R_0$-module for all $d \in \Z$.
In such situations, one can ask whether the $i$-th
cohomological Hilbert function $h_M^i$ of $M$ is of reverse
polynomial type: it always is when $R_0$ is Artinian, by
\cite[Theroem 17.1.9]{BS}. In this and the next section, we shall
study such situations where $R_0$ is not Artinian, and we shall
present an example in which $h_M^i$ is of reverse polynomial type,
and examples (in which $R_0$ is a polynomial ring over a field) in
which $h_M^i$ is not of reverse polynomial type: in this section
we give an example over each field $\Z/p\Z$, where $p$ is a prime
number; in the next section, we give an example over an arbitrary
field of characteristic zero.

In this section we study the cohomological Hilbert functions
arising from an example first studied by A.~Singh in
\cite{Singh00}; the example was also studied in \cite[\S2]{BKS},
where the asymptotic behaviour of the sets of associated primes of
the graded components of one of its local cohomology modules were
investigated.

In this section, we shall use several results from
\cite[\S2]{BKS}, and unexplained terminology will be found in that
section. Throughout this section, $L$ will denote either a field
or a principal ideal domain; let $R_0=L[X,Y,Z]$ and $R =
R_0[U,V,W]/(XU + YV + ZW)$, where $X,Y,Z,U,V,W$ are independent
indeterminates over $L$; we also assign degree $0$ to $X,Y,Z$ and
degree $1$ to $U,V,W$. Denote by $R_+$ the ideal of $R$ generated
by the images of $U, V, W$.

During the course of the section, we shall have occasion to take
$L$ to be $\Q$ and $\Z/p\Z$, where $p$ is a prime number. The way
in which Proposition 2.8 of \cite{BKS} is formulated means that
the calculations in that result can be used in these cases.

Notice that Theorem \ref{theorem1} implies that, when $L$ is a
field, for any $d\geq 3$, the vector space $H^3_{R+} (R)_{-d}$
over $L$ has finite dimension, and that $\dim_LH^3_{R+} (R)_{-d} =
\length _{R_0}H^3_{R_+}(M)_{-d}$.

\begin{defi} Suppose that $L$ is a field.
For all $d \in \Z$ we define $h_L(d)=\dim_{L} H^3_{R+}(R )_{d}.$
We abbreviate $h_{\Q}$ by $h_0$, and $h_{\Z/p\Z}$, for a prime number
$p$, by $h_p.$ Thus $h_0$ and $h_p$ are cohomological Hilbert functions.
\end{defi}

\begin{lem}\label{ch.2} Let $d \in \N$ with $d \geq 3$. Consider the matrix
$$T_d :=
\left[ \begin{array}{cccccc}
  A_{d-2} & XI_{d-2} & 0 & \dots &  & 0 \\
  0 & A_{d-3} & XI_{d-3} & 0 & \dots & 0 \\
  0 & 0 & A_{d-4} & XI_{d-4} &  & 0 \\
  \vdots & \vdots &  & \ddots & \ddots & \vdots \\
  0 & 0 & \dots & 0 & A_1 & XI_1
\end{array} \right]
$$
of\/ {\rm \cite[Lemma 2.2]{BKS}}. The induced monomorphism
$$
L[Y,Z]^{d-1 \choose 2}/\left(\Ima T_d \cap L[Y,Z]^{d-1 \choose
2}\right) \longrightarrow \Coker T_d
$$
is an $L$-isomorphism. Consequently, there is an $L$-isomorphism
$H^3_{R+} (R)_{-d} \cong \Coker H_d$, where $$H_d :=
\left[\begin{array}{ccccc}
A_{d-2} & 0                & 0 &  \ldots         &  0  \\
0       & A_{d-3} A_{d-2}  & 0 &  \ldots         &  0\\
0 & 0& A_{d-4}A_{d-3} A_{d-2} & \ldots           & 0\\
\vdots  &   \vdots    &             & \ddots    & \vdots \\

0        &  0    &  \ldots           &           & A_{1} A_{2}
\dots A_{d-2}
\end{array}\right],$$
as in\/ {\rm \cite[Theorem 2.4(ii)]{BKS}}.
\end{lem}

\begin{proof} The first part follows from the fact that, for all
integers $i$ with $1 \leq i \leq {d-1 \choose 2}$, the image of
$T_d$ contains a vector $X\mathbf{e}_i + \sum_{i<j\leq {d-1
\choose 2}} \alpha_j \mathbf{e}_j$, where each $\alpha_j$ is in
$L[Y,Z]$; we can use these to show that any member of
$L[X,Y,Z]^{d-1 \choose 2}$ differs by an element of $\Ima T_d$
from some element of $L[Y,Z]^{d-1 \choose 2}$.

The second part follows from the facts that the columns of $H_d$
generate $\Ima T_d \cap L[Y,Z]^{{d-1} \choose {2}}$ (by
\cite[Theorem 2.4(ii)]{BKS}) and that, by \cite[Lemma 2.2]{BKS},
there is an $L[X,Y,Z]$-isomorphism $\Coker T_d \cong H^3_{R+}(R
)_{-d}$.
\end{proof}

\begin{prop}\label{lemma4}\hfill
\begin{enumerate}
\item[(i)]
For any prime number $p$ and for all $d\geq 3$, we have $h_0(-d)
\leq h_p(-d)$; equality occurs if and only if
$$p\notin \Pi(d-2):=\left\{ p : p \mbox{~is a prime factor of~} {{d-2}
\choose {i}} \mbox{~for some~} i \in \{1, \dots, d-2\}\right\}.
$$
\item[(ii)] We have
$h_0(-d)=d(d-1)^2(d-2)/12$ for all $d\geq 3$. Thus the
cohomological Hilbert function $h_0$ is of reverse polynomial
type.
\end{enumerate}
\end{prop}
\begin{proof} We begin by considering the situation where $L$ is an
arbitrary field. As in \cite[Proposition 2.8]{BKS}, we consider
$L[Y,Z]$ as an $\nn^2$-graded ring in which $L[Y,Z]_{(0,0)} = L$
and $\deg Y^iZ^j = (i+j,i)$. We also endow the free
$L[Y,Z]$-module
$$
L[Y,Z]^r = L[Y,Z]\mathbf{e}_1 \oplus \cdots \oplus
L[Y,Z]\mathbf{e}_r
$$
with a structure as $\nn^2$-graded module over the $\nn^2$-graded
ring $L[Y,Z]$ in such a way that $\deg \mathbf{e}_i = (0,i)$ for
$i =1, \ldots, r$. Except where specified, we shall use this
bigrading henceforth in this proof.

Lemma \ref{ch.2} now implies that $h_L(-d)=\dim_{L} \Coker
H_d$ for all $d \geq 3$; hence, in the notation of
\cite[Proposition 2.8]{BKS}, for such a $d$ we have
$$h_L(-d)=\sum_{r=1}^{d-2} \dim_{L} \Coker Q_{r,d-1}.$$
We now use \cite[Proposition 2.8(ii)]{BKS} to calculate
$g_L(r,r+k) := \dim_L\Coker Q_{r,r+k}$ for $r,k \in \N$.

It follows from \cite[Proposition 2.8]{BKS} that
$\left(L[Y,Z]^r\right) _{(i,j)} = 0$ whenever $j > i+r$, that
$\Ima Q_{r,r+k}$ is a graded submodule of $L[Y,Z]^r$ and that
$(\Ima Q_{r,r+k})_{(i,j)} = \left(L[Y,Z]^r\right) _{(i,j)}$
whenever $i \geq 2k+r$. Hence
$$
\dim_L\Coker Q_{r,r+k} = \sum_{i=0}^{2k+r-1}\sum_{j \in \N}\dim_L
\left(L[Y,Z]^r\right) _{(i,j)} -
\sum_{i=0}^{2k+r-1}\sum_{j \in \N}\dim_L(\Ima Q_{r,r+k})_{(i,j)}.
$$
The first term on the right-hand side is equal to $r$ times the number of
monomials $Y^{\tau}Z^{\omega}$ (where $\tau, \omega \in \nn$) of
total (ordinary) degree not exceeding $2k +r -1$. Therefore
$$
\sum_{i=0}^{2k+r-1}\sum_{j \in \N}\dim_L
\left(L[Y,Z]^r\right) _{(i,j)}
= r\sum_{i=0}^{2k+r-1} (i+1) = \frac{r(2k+r)(2k+r+1)}{2}.
$$

Furthermore, it also follows from \cite[Proposition 2.8]{BKS} that
$(\Ima Q_{r,r+k})_{(i,j)} = 0$ whenever $i < k$, while if $k \leq i \leq
2k+r-1$, then $\dim_L(\Ima Q_{r,r+k})_{(i,j)}$ is equal to the rank
of a submatrix of $\widetilde{Q}_{r,r+k}$ (defined in
\cite[Proposition 2.8]{BKS}) made up
of the (consecutive) columns of that matrix numbered
$$\max\{j+k-i,1\},\max\{j+k-i,1\}+1, \ldots, \min \{j,r+k\}.$$
This rank depends on the characteristic of $L$: when $L = \Q$, each such
submatrix has maximal rank (by \cite[Corollary 2.12]{BKS}), whereas when
$L = \Z/p\Z$ for a prime number $p$, such a submatrix may not have
maximal rank. Therefore $g_{\Q}(r,r+k) \leq g_{\Z/p\Z}(r,r+k)$, and
equality holds if and only if, over $\Z/p\Z$, all submatrices of
$\widetilde{Q}_{r,r+k}$ formed by consecutive columns have maximal rank.
Let $p$ denote a prime number.
Now \cite[Corollary 2.14]{BKS} implies that when $p\notin \Pi(r+k-1)$
the above submatrices do all have maximal rank. Since
$$
h_L(-d) = \sum_{r=1}^{d-2} \dim_{L} \Coker Q_{r,d-1} \quad
\mbox{~for all~} d \geq 3,
$$
it follows that $h_0(-d) \leq h_p(-d)$ for all $d \geq 3$, and
that equality holds if $p\notin \Pi(d-2)$. However, when $p\in
\Pi(d-2)$, some of the $1 \times 1$ submatrices of
$$\widetilde{Q}_{1,d-1}=
\left[  {d-2 \choose 0} , {d-2 \choose 1} , \dots , {d-2 \choose \gamma} ,
\dots , {d-2  \choose 1},  {d-2  \choose 0}\right]$$
will vanish, and for such $p$ and $d$ we must have $h_0(d) < h_p(d)$.

(ii) Next, in the case where $L = \Q$, we calculate
$\sum_{i=0}^{2k+r-1}\sum_{j \in \N}\dim_{\Q}(\Ima Q_{r,r+k})_{(i,j)}.$
The comments in the preceding paragraph show that
$$\sum_{i=0}^{2k+r-1}\sum_{j \in \N}\dim_{\Q}(\Ima Q_{r,r+k})_{(i,j)} =
\sum_{i=k}^{2k+r-1}\sum_{j \in \N}\dim_{\Q}(\Ima
Q_{r,r+k})_{(i,j)},$$ and that, for an integer $i$ with $k \leq i
\leq 2k + r - 1$, the sum $\sum_{j \in \N}\dim_{\Q}(\Ima
Q_{r,r+k})_{(i,j)} =: T_i$ is equal to the sum of the ranks of the
submatrices of $Q_{r,r+k}$ obtained by selecting consecutive
columns numbered by the sets of integers in the list
$$\{1\}, \{1,2\}, \dots, \{1,2,\dots ,c\}, \{2,3,\dots, c+1\}, \dots,
\{r+k-c+1,r+k-c+2,\dots, r+k\},$$
$$ \{r+k-c+2,r+k-c+3,\dots, r+k\}, \dots,
\{r+k-1, r+k\} \textrm{\ and\ } \{r+k\}, $$
where $c=i-k+1$. Since each such
submatrix has maximal rank (by \cite[Corollary 2.12]{BKS}) and has $r$ rows,
it follows that
\begin{align*}
\sum_{i=k}^{2k+r-1}T_i & = \sum_{c=1}^r (k+r-c+1)c +
\sum_{c=r+1}^{k+r} (k+r-c+1)r + 2\sum_{c=1}^r \sum_{l=1}^{c-1} l +
2\sum_{c=r+1}^{k+r} \left( \sum_{l=1}^{r} l + \sum_{l=r+1}^{c-1} r \right) \\
& = \frac{(k+r)^2(k+r+1)}{2} - \frac{k^2(k+1)}{2}.
\end{align*}
Hence
\begin{align*}
g_{\Q}(r,r+k) & = \dim_{\Q}\Coker Q_{r,r+k} \\
& =
\sum_{i=0}^{2k+r-1}\sum_{j \in \N}\dim_{\Q}
\left(L[Y,Z]^r\right) _{(i,j)} -
\sum_{i=0}^{2k+r-1}\sum_{j \in \N}\dim_{\Q}(\Ima Q_{r,r+k})_{(i,j)} \\
& =
\frac{r(2k+r)(2k+r+1)}{2} - \frac{(k+r)^2(k+r+1)}{2} + \frac{k^2(k+1)}{2} \\
& = \frac{rk(k+r)}{2}.
\end{align*}

Finally, we return to our cohomological Hilbert function: for all
$d \geq 3$,
\begin{align*}
h_0(-d) & = \sum_{r=1}^{d-2} \dim_{\Q} \Coker Q_{r,d-1} \\
& = \sum_{r=1}^{d-2} \frac{r(d-1-r)(d-1)}{2}\\
& =\frac{d(d-1)^2(d-2)}{12}.
\end{align*}
\end{proof}

\begin{thm}\label{theorem2}\hfill
\begin{itemize}
\item[(i)] For any prime number $p$ both the sets
$$\{ d : d\geq 3,~ h_p(-d)=h_0(-d)\}\quad \mbox{~and~} \quad
\{ d : d\geq 3,~ h_p(-d)> h_0(-d)\}$$ are infinite.

\item[(ii)] None of the cohomological Hilbert functions
$h_p$ is of reverse polynomial type for any prime number $p$.

\item[(iii)] If $p$ and $q$ are different prime numbers then the set
$\{ d : d\geq 3,~ h_p(-d)> h_q(-d)\}$ is infinite.
\end{itemize}
\end{thm}
\begin{proof} (i)
This follows from Proposition \ref{lemma4}(i) and \cite[Lemma
2.16]{BKS}.

(ii) Let $p$ be any prime number. If $h_p(r)=P(r)$ for some
polynomial $P\in \Q[T]$ and for all $r < < -3$, then $P(r)=h_0(r)$
for infinitely many $r < < -3$ by part (i). Thus, by Proposition
\ref{lemma4}(ii), we must have $P = T(T-1)^2(T-2)/12$, so that the
set $\{ d : d\geq 3,~ h_p(-d)> h_0(-d)\}$ is finite. But this
contradicts part (i).

(iii) Assume now that $p$ and $q$ are different prime numbers. By
Proposition \ref{lemma4}(i), it is enough to show that $\{ j\in \N
: p\in \Pi(j),\ q\notin \Pi(j) \}$ is infinite. The proof of Lemma
2.16 in \cite{BKS} shows that $q\notin \Pi(q^k-1)$ for all $k\geq
1$. On the other hand, if $p$ divides $a\in \N$ then $p$ divides
${a \choose 1}=a$, so that $p\in \Pi(a)$; it is therefore enough
to show that $p$ divides $q^k-1$ for infinitely many $k\in \N$.
Let $\alpha$ be the order of $q$ in the multiplicative group of
$\Z/p\Z$. For all $\beta\in \N$ we have
$$q^{\alpha p^\beta} -1 \equiv q^\alpha -1 \equiv 0\mod p.$$
\end{proof}

\section{\bf An example in characteristic zero}

In \S\ref{sec2}, we provided, for each prime number $p$, an
example of a cohomological Hilbert function of a standard
positively graded finitely generated algebra over the field
$\Z/p\Z$ that fails to be of reverse polynomial type. In this
section, we provide an example over a field of characteristic $0$
that exhibits similar behaviour.

Fix $K$ to be any field of characteristic zero. Let $R_0=K[X,Y]$
and let $S=R_0[U,V]$, where $X,Y,U,V$ are independent
indeterminates over $K$. Define a grading on $S$ by declaring that
$\deg X=\deg Y=0$ and $\deg U=\deg V=1$. Let $f=2X^2V^2 + 2 XY UV
+ Y^2 U^2$ and let $R=S/fS$. Notice that $f$ is homogeneous and
hence $R$ is graded. Let $S_+$ be the ideal of $S$ generated by
$U$ and $V$ and let $R_+$ be the ideal of $R$ generated by the
images of $U$ and $V$.

We will study the graded components of $H^2_{R_+}(R)$ by
exploiting the fact that this local cohomology module is
homogeneously isomorphic to $H^2_{S_+}(S)/fH^2_{S_+}(S)$, and that
$H^2_{S_+}(S)$ can be realised as the module $R_0[U^-,V^-]$ of
inverse polynomials. Thus, for all $d \in \Z$, the $(-d)$-th
graded component $H^2_{R_+}(R)_{-d}$ of $H^2_{R_+}(R)$ is
isomorphic to the cokernel of the $R_0$-homomorphism
$$
f_{d}: R_0[U^-,V^-]_{-d-2} \lra R_0[U^-,V^-]_{-d}
$$
given by multiplication by $f$, as described in Section
\ref{section1}. Note that $H^2_{R_+}(R)_{r} = 0$ for all $r >-2$,
and that Theorem \ref{theorem1} shows that, for all $d$, the
$R_0$-module $H^2_{R_+}(R)_{-d}$ has finite length
$$
h^2_R(-d) := \length_{R_0}H^2_{R_+}(R)_{-d} =
\dim_KH^2_{R_+}(R)_{-d}\mbox{;}
$$
Thus $h^2_R : \Z \rightarrow \nn$ is a cohomological Hilbert function,
as explained at the beginning of the last section. Our aim in this
section is to show that $h^2_R$ is not of reverse polynomial type.

If we use the ordering of bases described in Section
\ref{section1} (with $U > V$) for both the source and target of
each $f_d$, we can see that each $f_d~(d\geq 2)$ is given by
multiplication on the left by the $(d-1) \times (d+1)$ tridiagonal
matrix
$$
C_{d-1}:= \left[
\begin{array}{cccccccc}
  2X^2 & 2 XY & Y^2 & 0 & 0 & \dots  & 0 & 0 \\
  0 & 2X^2 & 2XY & Y^2 & 0 & \dots & 0 & 0 \\
  0 & 0 & 2X^2 & 2XY & Y^2 & \dots & 0 & 0 \\
  \vdots  & \vdots  & & \ddots & \ddots & & \vdots& \vdots \\
   \vdots & \vdots  & & & \ddots & \ddots & \vdots& \vdots  \\
  0 & 0 & \dots & 0 & 0 & 2X^2 & 2XY & Y^2  \\
\end{array}
\right] .
$$
We also write
$$
\widetilde{C}_{d-1}:= \left[
\begin{array}{cccccccc}
  2 & 2 & 1 & 0 & 0 & \dots & 0 & 0 \\
  0 & 2 & 2 & 1 & 0 & \dots & 0 & 0 \\
  0 & 0 & 2 & 2 & 1 & \dots & 0 & 0 \\
  \vdots  & \vdots  & & \ddots & \ddots & &\vdots & \vdots \\
   \vdots & \vdots  & & & \ddots & \ddots & \vdots& \vdots  \\
  0 & 0 & \dots & 0 & 0 & 2 & 2 & 1  \\
\end{array}
\right] ,
$$
the result of evaluation of $C_{d-1}$ at $X=Y=1$.

\begin{lem}\label{lemma5}\hfil
\begin{enumerate}
\item[(i)]
For each $n \in \N$, let $\widetilde{D}_n$ be the $n\times n$ tridiagonal
matrix
$$\left[
\begin{array}{cccccc}
   2 & 1 & 0 & 0 & \dots & 0  \\
   2 & 2 & 1 & 0 & \dots & 0  \\
   0 & 2 & 2 & 1 & \dots & 0  \\
   \vdots   & \vdots & \ddots & \ddots & &\vdots  \\
   \vdots   &\vdots  & & \ddots & \ddots & \vdots  \\
   0 & \dots & 0 & 0 & 2 & 2   \\
\end{array}
\right]
$$
obtained by removing the first and last columns of
$\widetilde{C}_{n}$. Then $\det \widetilde{D}_n=0$ if and only if
$n\equiv 3 \mod 4$.

\item[(ii)]
Any submatrix of $\widetilde{C}_{n}$ consisting of consecutive
columns has maximal rank\/ \emph{except} for the matrix obtained from
$\widetilde{C}_n$ by removing its first and last
columns when $n\equiv 3 \mod 4$.
This exceptional submatrix has rank $n-1$.
\end{enumerate}
\end{lem}
\begin{proof} (i)
Write $\Delta_n=\det \widetilde{D}_n$ for all $n\geq 1$. We have
$\Delta_1= 2$ and $\Delta_2=2$; for $n\geq 3$ we can expand the
determinant $\Delta_n$ by the first row of
$\widetilde{D}_n$ to obtain $\Delta_n=2
\Delta_{n-1} - 2 \Delta_{n-2}$. So for all $n\geq 3$ we can write
$$
\left[ \begin{array}{c} \Delta_n\\ \Delta_{n-1} \end{array} \right]=
\left[ \begin{array}{cc} 2&-2\\1&0 \end{array} \right]
\left[ \begin{array}{c} \Delta_{n-1}\\ \Delta_{n-2} \end{array} \right]
$$
and by induction
$$
\left[ \begin{array}{c} \Delta_n\\ \Delta_{n-1} \end{array} \right]=
\left[ \begin{array}{cc} 2&-2\\1&0 \end{array} \right]^{n-2}
\left[ \begin{array}{c} 2\\ 2 \end{array} \right] .
$$
The $2 \times 2$ rational matrix
$\left[ \begin{array}{cc} 2&-2\\1&0 \end{array} \right]$
has complex eigenvalues
$\sqrt{2} e^{i \pi/4}$ and $\sqrt{2} e^{-i \pi/4}$, and by
diagonalizing it we see that, for $n \geq 3$,
\begin{align*}
\left[ \begin{array}{c} \Delta_n\\ \Delta_{n-1} \end{array} \right] & =
\left[ \begin{array}{cc} 2&-2\\1&0 \end{array} \right]^{n-2}
\left[ \begin{array}{c} 2\\ 2 \end{array} \right]\\
& =  (\sqrt{2})^{n-2}
\left[ \begin{array}{cc} 1+i&1-i\\1&1 \end{array} \right]
\left[ \begin{array}{cc} e^{i (n-2) \pi/4}&0\\ 0&e^{-i (n-2) \pi/4}
\end{array} \right]
\left[ \begin{array}{cc} -i&1+i\\i&1-i \end{array} \right]
\left[ \begin{array}{c} 1\\ 1 \end{array} \right]\\
& = (\sqrt{2})^{n-2}
\left[ \begin{array}{cc} 1+i&1-i\\1&1 \end{array} \right]
\left[ \begin{array}{c} e^{i (n-2) \pi/4}\\ e^{-i (n-2) \pi/4} \end{array}
\right]\\
& = (\sqrt{2})^{n-2} \left[ \begin{array}{c}
(1+i) e^{i (n-2) \pi/4}+ (1-i) e^{-i (n-2) \pi/4} \\
e^{i (n-2) \pi/4} + e^{-i (n-2) \pi/4}
\end{array} \right].
\end{align*}
Hence $\Delta _n = 0$ if and only if $n\equiv 3\mod 4$.

(ii) The matrix obtained from $\widetilde{C}_{n}$ by removal of
its first and last columns is $\widetilde{D}_n$, and it is
straightforward to check that all other selections of $n$ or fewer
consecutive columns from $\widetilde{C}_{n}$ are linearly
independent.
\end{proof}

\begin{thm}\label{theorem3}
The cohomological Hilbert function $h^2_R:\Z \rightarrow \nn$ is such
that, for all $d\geq 2$,
$$h^2_R(-d)=
\begin{cases}
d^2-1 & \text{if $d \not \equiv 0\mod 4$}, \\
d^2 & \text{if $d \equiv 0 \mod 4$}.
\end{cases}
$$
Hence $h^2_R$ is not of reverse polynomial type.
\end{thm}
\begin{proof}
Fix an integer $d \geq 2$. The length $h^2_R(-d)$ is $\dim _K\Coker
C_{d-1}$. In order to calculate this dimension, we consider
$K[X,Y]$ as an $\nn^2$-graded ring in which $K[X,Y]_{(0,0)} = K$
and $\deg X^iY^j = (i+j,j)$ for all $(i,j) \in \nn^2$. For
convenience, we set $d-1 =: n$. Turn the free $R_0$-module
$$
R_0^n = K[X,Y]\mathbf{e}_1 \oplus \cdots \oplus
K[X,Y]\mathbf{e}_n
$$
into an $\nn^2$-graded module over the $\nn^2$-graded ring
$R_0 = K[X,Y]$ in such a way that $\deg \mathbf{e}_i = (0,i)$ for $i =1,
\ldots, n$. All references to gradings in the rest of this
proof refer to this $\nn^2$-grading. Note that $(R_0^n)_{(i,j)} = 0$
whenever $j > n+i$.

For each $j = 1, \ldots, n+2$, let $\mathbf{c}_j$ denote the $j$-th column of
$C_n$, and note that $\mathbf{c}_j$ is homogeneous of degree $(2,j)$. Thus
$\Ima C_n$, the $R_0$-submodule of $R_0^n$ generated by the columns of $C_n$,
is graded; hence $\Coker C_n$ is graded, too. Note that
$(\Ima C_n)_{(i,j)} = 0$ whenever $i < 2$.

It is not hard to see that the vectors
$$X^2 \mathbf{e}_1, X^3 \mathbf{e}_2, \dots, X^{n+1} \mathbf{e}_{n}
\quad \textrm{\ and\ } \quad
Y^2 \mathbf{e}_{n}, Y^3 \mathbf{e}_{n-1}, \dots, Y^{n+1} \mathbf{e}_{1} $$
are in the image of $C_n$.
Hence $\left(\Ima C_{n}\right)_{(i,j)}=\left(R_0^n\right)_{(i,j)}$
whenever $i>n+1$.

These observations lead to the conclusion that
$$
\dim_K\Coker C_{n} = \sum_{i=0}^{n+1}\sum_{j \in \N}\dim_K
\left(K[X,Y]^n\right) _{(i,j)} -
\sum_{i=2}^{n+1}\sum_{j \in \N}\dim_K(\Ima C_{n})_{(i,j)}.
$$
The first term on the right-hand side is equal to $n$ times the number of
monomials $X^{\tau}Y^{\omega}$ (where $\tau, \omega \in \nn$) of
total (ordinary) degree not exceeding $n+1$. Therefore
$$
\sum_{i=0}^{n+1}\sum_{j \in \N}\dim_K \left(K[X,Y]^n\right)
_{(i,j)} = n\sum_{i=0}^{n+1} (i+1) = \frac{n(n+2)(n+3)}{2}.
$$

Now choose an integer $i$ with $2 \leq i \leq n+1$. Our observations above
show that, for any $j \in \N$,
$$\left(\Ima C_n\right)_{(i,j)}= \sum_{\sigma =
\max\{2+j-i,1\}}^{\min\{j,n+2\}}
KX^{i-2-j+\sigma}Y^{j-\sigma}\mathbf{c}_{\sigma}.$$ Thus
$\dim_K(\Ima C_n)_{(i,j)}$ is equal to the rank of the submatrix
of $\widetilde{C}_{n}$ made up of the (consecutive) columns of
that matrix numbered
$$\max\{2+j-i,1\},\max\{2+j-i,1\}+1, \ldots, \min \{j,n+2\}.$$
It follows that $\sum_{j \in \N}\dim_{K}(\Ima C_{n})_{(i,j)} $ is
equal to the sum of the ranks of the submatrices of
$\widetilde{C}_{n}$ obtained by selecting consecutive columns
numbered by the sets of integers in the list
$$\{1\}, \{1,2\}, \dots, \{1,2,\dots ,i-1\}, \{2,3,\dots, i\}, \dots,
\{n-i+4,n-i+5,\dots, n+2\},$$
$$ \{n-i+5,n-i+6,\dots, n+2\}, \dots,
\{n+1, n+2\} \textrm{\ and\ } \{n+2\}. $$
By Lemma \ref{lemma5}(ii), all these submatrices have maximal rank, except
when $n\equiv 3 \mod 4$ and $i = n+1$,
when they all have maximal rank except for
that corresponding to the choice
$\{2,3,\dots, n+1\}$, which has rank $n-1$ rather than $n$.

It follows that, unless $n\equiv 3 \mod 4$ and $i = n+1$,
\begin{align*}
\sum_{j \in \N}\dim_{K}(\Ima C_{n})_{(i,j)} & =
2\sum_{\alpha = 1}^{i-2} \alpha + (n-i+4)(i-1) \\
& = (i-2)(i-1) + (n-i+4)(i-1) = (i-1)(n+2)\mbox{;}
\end{align*}
in the exceptional case, the sum is one less than that given by
the above formula. Hence, unless $n\equiv 3 \mod 4$,
\begin{align*}
\dim_K\Coker C_{n} & = \sum_{i=0}^{n+1}\sum_{j \in \N}\dim_K
\left(K[X,Y]^n\right) _{(i,j)} -
\sum_{i=2}^{n+1}\sum_{j \in \N}\dim_K(\Ima C_{n})_{(i,j)}\\
& = \frac{n(n+2)(n+3)}{2} - \sum_{i=2}^{n+1} (i-1)(n+2)\\
& = \frac{n(n+2)(n+3)}{2} - \frac{n(n+1)(n+2)}{2} = n(n+2)\mbox{;}
\end{align*}
when $n\equiv 3 \mod 4$, we have $\dim_K\Coker C_{n} = n(n+2) +
1$.

Thus, since $n = d-1$, we have shown that
$$
h^2_R(-d) = \dim _K\Coker C_{d-1} = \begin{cases}
d^2-1 & \text{if $d \not \equiv 0 \mod 4$}, \\
d^2 & \text{if $d \equiv 0 \mod 4$}.
\end{cases}$$

If $h^2_R(r)$ were to agree with $P(r)$ for a polynomial $P \in \Q[T]$
for all $r < < 0$, then $P$ would have to be both $T^2-1$ and
$T^2$, which is, of course, absurd.
\end{proof}

\section*{\bf Acknowledgement}
We would like to thank Markus Brodmann for helpful discussions about the
content of this paper.

\end{document}